\documentclass[11pt]{article}
\usepackage[utf8]{inputenc}
\usepackage{amsmath, amssymb, amsthm}
\usepackage{geometry}
\usepackage{tikz}
\usetikzlibrary{trees, positioning, arrows.meta, calc, decorations.pathreplacing}
\usepackage{hyperref}
\usepackage{enumitem}

\newtheorem{theorem}{Theorem}
\newtheorem{lemma}[theorem]{Lemma}
\newtheorem{proposition}[theorem]{Proposition}

\newtheorem{conjecture}[theorem]{Conjecture}
\theoremstyle{definition}
\newtheorem{definition}[theorem]{Definition}

\theoremstyle{remark}

\tikzset{
    treenode/.style = {circle, fill=black, inner sep=0pt, minimum size=5pt},
    treeedge/.style = {thick},
    node/.style = {circle, fill=black, inner sep=2pt},
}

\title{\textbf{A Lower Bound for Kruskal's Weak Tree Function:} \\ \textbf{$\mathrm{tree}(3) \geq 844{,}424{,}930{,}131{,}960$}}
\author{Mark Giroux}
\date{}

\begin{document}

\maketitle

\begin{abstract}
We establish an explicit lower bound for Kruskal's weak tree function at $n=3$, proving that $\mathrm{tree}(3) \geq 844{,}424{,}930{,}131{,}960 = 3 \cdot 2^{48} - 8$. This is achieved by constructing an explicit sequence of unlabeled rooted trees satisfying the constraints of the weak tree function and carefully analyzing the combinatorics of the ``leg elimination'' process. Our bound significantly exceeds previous estimates and demonstrates that even for small arguments, the weak tree function exhibits rapid growth.
\end{abstract}

\section{Introduction}

Kruskal's tree theorem \cite{kruskal1960} is a landmark result in combinatorics and proof theory, stating that finite trees over a well-quasi-ordered label set are themselves well-quasi-ordered under homeomorphic embedding. A consequence of this theorem is that certain sequences of trees must eventually contain an embedding, which led Friedman to define fast-growing functions measuring the longest possible ``bad'' sequences.

The \textit{weak tree function} $\mathrm{tree}(n)$ considers unlabeled trees (equivalently, trees with a single label) and asks for the longest sequence avoiding homeomorphic embeddings. While the labeled TREE function has received considerable attention---particularly TREE(3), which vastly exceeds Graham's number---the weak tree function has been comparatively less studied.

Previous estimates for $\mathrm{tree}(3)$ have varied widely. Friedman conjectured that the related function $\mathrm{FFF}(2) = \mathrm{tree}(3)$ would be less than 100 \cite{friedman2006}. In this paper, we prove that the true value is at least $844{,}424{,}930{,}131{,}960$, a number exceeding $8 \times 10^{14}$.

\section{Preliminaries}

\begin{definition}[Rooted Tree]
A \textit{rooted tree} is a connected acyclic graph $T = (V, E)$ with a distinguished vertex $r \in V$ called the \textit{root}. For vertices $u, v \in V$, we say $v$ is a \textit{descendant} of $u$ if $u$ lies on the unique path from $r$ to $v$. We say $v$ is a \textit{child} of $u$ if $v$ is a descendant of $u$ and $(u, v) \in E$.
\end{definition}

\begin{definition}[Infimum / Least Common Ancestor]
For vertices $u, v$ in a rooted tree $T$, the \textit{infimum} (or \textit{least common ancestor}) $\inf(u, v)$ is the unique vertex $w$ such that $w$ is an ancestor of both $u$ and $v$, and no descendant of $w$ has this property.
\end{definition}

\begin{definition}[Inf-Embedding]
Let $T_1 = (V_1, E_1, r_1)$ and $T_2 = (V_2, E_2, r_2)$ be rooted trees. An \textit{inf-embedding} (or \textit{homeomorphic embedding preserving infima}) is an injective map $f: V_1 \to V_2$ such that:
\begin{enumerate}[label=(\roman*)]
    \item If $v$ is a descendant of $u$ in $T_1$, then $f(v)$ is a descendant of $f(u)$ in $T_2$.
    \item For all $u, v \in V_1$, we have $f(\inf_{T_1}(u, v)) = \inf_{T_2}(f(u), f(v))$.
\end{enumerate}
If such an embedding exists, we write $T_1 \leq T_2$.
\end{definition}

\begin{definition}[Weak Tree Function]
The \textit{weak tree function} $\mathrm{tree}(n)$ is defined as the length of the longest sequence $(T_1, T_2, \ldots, T_m)$ of unlabeled rooted trees such that:
\begin{enumerate}[label=(\roman*)]
    \item For each $k \geq 1$, the tree $T_k$ has at most $k + n$ vertices.
    \item For all $i < j$, we have $T_i \not\leq T_j$ (no tree inf-embeds into a later tree).
\end{enumerate}
\end{definition}

\subsection{Tree Representation}

We represent rooted trees using visual diagrams where vertices are shown as black dots connected by edges. This notation makes the structure immediately clear and facilitates understanding of homeomorphic embeddings.

\begin{definition}[Visual Tree Representation]
A rooted tree is represented as:
\begin{itemize}
    \item Vertices are drawn as filled black circles (nodes).
    \item Edges are drawn as thick lines connecting parent nodes to their children.
    \item The root node is positioned at the top, with children below.
    \item For trees with long descending chains, we use dotted lines with a number indicating the chain length to save space.
\end{itemize}
\end{definition}

\section{Known Values}

The following values are easily verified:

\begin{proposition}
$\mathrm{tree}(1) = 2$, achieved by the sequence:
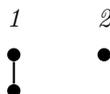
\begin{center}
\begin{tikzpicture}[scale=0.8]
    % Tree 1
    \begin{scope}[shift={(0,0)}]
        \node[above] at (0,0.3) {\small 1};
        \node[treenode] (r1) at (0,0) {};
        \node[treenode] (a1) at (0,-0.6) {};
        \draw[treeedge] (r1) -- (a1);
    \end{scope}

    % Tree 2
    \begin{scope}[shift={(1.5,0)}]
        \node[above] at (0,0.3) {\small 2};
        \node[treenode] (r2) at (0,0) {};
    \end{scope}
\end{tikzpicture}
\end{center}
A root with one child (2 vertices), followed by a single root node (1 vertex).
\end{proposition}

\begin{proposition}
$\mathrm{tree}(2) = 5$, achieved by the sequence:
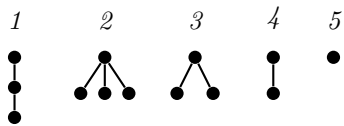
\begin{center}
\begin{tikzpicture}[scale=0.8]
    % Tree 1: chain of 3
    \begin{scope}[shift={(0,0)}]
        \node[above] at (0,0.3) {\small 1};
        \node[treenode] (r1) at (0,0) {};
        \node[treenode] (a1) at (0,-0.5) {};
        \node[treenode] (b1) at (0,-1.0) {};
        \draw[treeedge] (r1) -- (a1);
        \draw[treeedge] (a1) -- (b1);
    \end{scope}

    % Tree 2: root with 3 leaves
    \begin{scope}[shift={(1.5,0)}]
        \node[above] at (0,0.3) {\small 2};
        \node[treenode] (r2) at (0,0) {};
        \node[treenode] (a2) at (-0.4,-0.6) {};
        \node[treenode] (b2) at (0,-0.6) {};
        \node[treenode] (c2) at (0.4,-0.6) {};
        \draw[treeedge] (r2) -- (a2);
        \draw[treeedge] (r2) -- (b2);
        \draw[treeedge] (r2) -- (c2);
    \end{scope}

    % Tree 3: root with 2 leaves
    \begin{scope}[shift={(3.0,0)}]
        \node[above] at (0,0.3) {\small 3};
        \node[treenode] (r3) at (0,0) {};
        \node[treenode] (a3) at (-0.3,-0.6) {};
        \node[treenode] (b3) at (0.3,-0.6) {};
        \draw[treeedge] (r3) -- (a3);
        \draw[treeedge] (r3) -- (b3);
    \end{scope}

    % Tree 4: root with 1 leaf
    \begin{scope}[shift={(4.3,0)}]
        \node[above] at (0,0.3) {\small 4};
        \node[treenode] (r4) at (0,0) {};
        \node[treenode] (a4) at (0,-0.6) {};
        \draw[treeedge] (r4) -- (a4);
    \end{scope}

    % Tree 5: single node
    \begin{scope}[shift={(5.3,0)}]
        \node[above] at (0,0.3) {\small 5};
        \node[treenode] (r5) at (0,0) {};
    \end{scope}
\end{tikzpicture}
\end{center}
Chain of 3 nodes, followed by root with 3, 2, then 1 leaf children, ending with a single node.
\end{proposition}

\section{Main Result}

\begin{theorem}\label{thm:main}
$\mathrm{tree}(3) \geq 844{,}424{,}930{,}131{,}960 = 3 \cdot 2^{48} - 8$.
\end{theorem}

The proof proceeds by constructing an explicit valid sequence and carefully counting its length.

\subsection{The Initial Sequence}

We begin at step 4 (since $\mathrm{tree}(3)$ allows $k+3$ vertices at step $k$, so step 1 allows 4 vertices). The first four trees are:

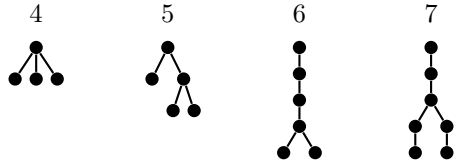
\begin{center}
\begin{tikzpicture}[scale=0.7]
    % Tree 4: root with 3 children
    \begin{scope}[shift={(0,0)}]
        \node[above] at (0,0.3) {\small 4};
        \node[treenode] (r4) at (0,0) {};
        \node[treenode] (a4) at (-0.4,-0.6) {};
        \node[treenode] (b4) at (0,-0.6) {};
        \node[treenode] (c4) at (0.4,-0.6) {};
        \draw[treeedge] (r4) -- (a4);
        \draw[treeedge] (r4) -- (b4);
        \draw[treeedge] (r4) -- (c4);
    \end{scope}

    % Tree 5: root with 2 children, right has 2 children
    \begin{scope}[shift={(2.5,0)}]
        \node[above] at (0,0.3) {\small 5};
        \node[treenode] (r5) at (0,0) {};
        \node[treenode] (a5) at (-0.3,-0.6) {};
        \node[treenode] (b5) at (0.3,-0.6) {};
        \node[treenode] (c5) at (0.1,-1.2) {};
        \node[treenode] (d5) at (0.5,-1.2) {};
        \draw[treeedge] (r5) -- (a5);
        \draw[treeedge] (r5) -- (b5);
        \draw[treeedge] (b5) -- (c5);
        \draw[treeedge] (b5) -- (d5);
    \end{scope}

    % Tree 6: 4 descending nodes, bottom has 2 children
    \begin{scope}[shift={(5,0)}]
        \node[above] at (0,0.3) {\small 6};
        \node[treenode] (r6) at (0,0) {};
        \node[treenode] (a6) at (0,-0.5) {};
        \node[treenode] (b6) at (0,-1.0) {};
        \node[treenode] (c6) at (0,-1.5) {};
        \node[treenode] (d6) at (-0.3,-2.0) {};
        \node[treenode] (e6) at (0.3,-2.0) {};
        \draw[treeedge] (r6) -- (a6);
        \draw[treeedge] (a6) -- (b6);
        \draw[treeedge] (b6) -- (c6);
        \draw[treeedge] (c6) -- (d6);
        \draw[treeedge] (c6) -- (e6);
    \end{scope}

    % Tree 7: two-leg structure
    \begin{scope}[shift={(7.5,0)}]
        \node[above] at (0,0.3) {\small 7};
        \node[treenode] (r7) at (0,0) {};
        \node[treenode] (a7) at (0,-0.5) {};
        \node[treenode] (b7) at (0,-1.0) {};
        \node[treenode] (c7) at (-0.3,-1.5) {};
        \node[treenode] (d7) at (0.3,-1.5) {};
        \node[treenode] (e7) at (-0.3,-2.0) {};
        \node[treenode] (f7) at (0.3,-2.0) {};
        \draw[treeedge] (r7) -- (a7);
        \draw[treeedge] (a7) -- (b7);
        \draw[treeedge] (b7) -- (c7);
        \draw[treeedge] (b7) -- (d7);
        \draw[treeedge] (c7) -- (e7);
        \draw[treeedge] (d7) -- (f7);
    \end{scope}
\end{tikzpicture}
\end{center}

Trees 8--11 show the systematic leg reduction:

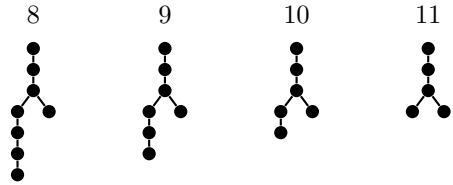
\begin{center}
\begin{tikzpicture}[scale=0.7]
    % Tree 8
    \begin{scope}[shift={(0,0)}]
        \node[above] at (0,0.3) {\small 8};
        \node[treenode] (r) at (0,0) {};
        \node[treenode] (a) at (0,-0.4) {};
        \node[treenode] (b) at (0,-0.8) {};
        \node[treenode] (c) at (-0.3,-1.2) {};
        \node[treenode] (d) at (0.3,-1.2) {};
        \node[treenode] (e) at (-0.3,-1.6) {};
        \node[treenode] (f) at (-0.3,-2.0) {};
        \node[treenode] (g) at (-0.3,-2.4) {};
        \draw[treeedge] (r) -- (a);
        \draw[treeedge] (a) -- (b);
        \draw[treeedge] (b) -- (c);
        \draw[treeedge] (b) -- (d);
        \draw[treeedge] (c) -- (e);
        \draw[treeedge] (e) -- (f);
        \draw[treeedge] (f) -- (g);
    \end{scope}

    % Tree 9
    \begin{scope}[shift={(2.5,0)}]
        \node[above] at (0,0.3) {\small 9};
        \node[treenode] (r) at (0,0) {};
        \node[treenode] (a) at (0,-0.4) {};
        \node[treenode] (b) at (0,-0.8) {};
        \node[treenode] (c) at (-0.3,-1.2) {};
        \node[treenode] (d) at (0.3,-1.2) {};
        \node[treenode] (e) at (-0.3,-1.6) {};
        \node[treenode] (f) at (-0.3,-2.0) {};
        \draw[treeedge] (r) -- (a);
        \draw[treeedge] (a) -- (b);
        \draw[treeedge] (b) -- (c);
        \draw[treeedge] (b) -- (d);
        \draw[treeedge] (c) -- (e);
        \draw[treeedge] (e) -- (f);
    \end{scope}

    % Tree 10
    \begin{scope}[shift={(5,0)}]
        \node[above] at (0,0.3) {\small 10};
        \node[treenode] (r) at (0,0) {};
        \node[treenode] (a) at (0,-0.4) {};
        \node[treenode] (b) at (0,-0.8) {};
        \node[treenode] (c) at (-0.3,-1.2) {};
        \node[treenode] (d) at (0.3,-1.2) {};
        \node[treenode] (e) at (-0.3,-1.6) {};
        \draw[treeedge] (r) -- (a);
        \draw[treeedge] (a) -- (b);
        \draw[treeedge] (b) -- (c);
        \draw[treeedge] (b) -- (d);
        \draw[treeedge] (c) -- (e);
    \end{scope}

    % Tree 11
    \begin{scope}[shift={(7.5,0)}]
        \node[above] at (0,0.3) {\small 11};
        \node[treenode] (r) at (0,0) {};
        \node[treenode] (a) at (0,-0.4) {};
        \node[treenode] (b) at (0,-0.8) {};
        \node[treenode] (c) at (-0.3,-1.2) {};
        \node[treenode] (d) at (0.3,-1.2) {};
        \draw[treeedge] (r) -- (a);
        \draw[treeedge] (a) -- (b);
        \draw[treeedge] (b) -- (c);
        \draw[treeedge] (b) -- (d);
    \end{scope}
\end{tikzpicture}
\end{center}

Tree 12 enters a symmetric two-leg configuration with depth-4 legs:

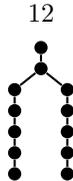
\begin{center}
\begin{tikzpicture}[scale=0.7]
    \node[above] at (0,0.3) {\small 12};
    \node[treenode] (r) at (0,0) {};
    \node[treenode] (a) at (0,-0.4) {};
    % Left leg
    \node[treenode] (l1) at (-0.5,-0.8) {};
    \node[treenode] (l2) at (-0.5,-1.2) {};
    \node[treenode] (l3) at (-0.5,-1.6) {};
    \node[treenode] (l4) at (-0.5,-2.0) {};
    \node[treenode] (l5) at (-0.5,-2.4) {};
    % Right leg
    \node[treenode] (r1) at (0.5,-0.8) {};
    \node[treenode] (r2) at (0.5,-1.2) {};
    \node[treenode] (r3) at (0.5,-1.6) {};
    \node[treenode] (r4) at (0.5,-2.0) {};
    \node[treenode] (r5) at (0.5,-2.4) {};
    \draw[treeedge] (r) -- (a);
    \draw[treeedge] (a) -- (l1);
    \draw[treeedge] (a) -- (r1);
    \draw[treeedge] (l1) -- (l2);
    \draw[treeedge] (l2) -- (l3);
    \draw[treeedge] (l3) -- (l4);
    \draw[treeedge] (l4) -- (l5);
    \draw[treeedge] (r1) -- (r2);
    \draw[treeedge] (r2) -- (r3);
    \draw[treeedge] (r3) -- (r4);
    \draw[treeedge] (r4) -- (r5);
\end{tikzpicture}
\end{center}

\subsection{The Leg Elimination Process}

The leg elimination process systematically reduces a symmetric two-leg configuration to simpler forms while avoiding embeddings. The process works as follows:

\begin{enumerate}
    \item Starting from a symmetric configuration $(x, x)$, we remove one node from the right leg, giving $(x, x-1)$.
    \item We then add the maximum number of nodes allowed by the vertex constraint to the left leg, creating an asymmetric configuration with a much longer left leg.
    \item At each subsequent step, we remove one node from the left leg until both legs again have matching depth, reaching the symmetric configuration $(x-1, x-1)$.
    \item The process repeats from step 1 with the new symmetric configuration, continuing until we reach the base case.
\end{enumerate}

This recursive structure leads to the following formula for the total number of steps:

\begin{lemma}[Leg Elimination Formula]\label{lem:leg}
The number of steps required to eliminate a symmetric two-leg configuration of depth $x$ is:
\[
L(x) = 6 \cdot 2^x - 2x - 6
\]
\end{lemma}

Applying Lemma~\ref{lem:leg} with $x = 4$:
\[
L(4) = 6 \cdot 2^4 - 2(4) - 6 = 96 - 8 - 6 = 82 \text{ steps}
\]

Adding the initial 12 steps: $12 + 82 = 94$.

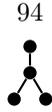
\begin{center}
\begin{tikzpicture}[scale=0.7]
    \node[above] at (0,0.3) {\small 94};
    \node[treenode] (r) at (0,0) {};
    \node[treenode] (a) at (0,-0.5) {};
    \node[treenode] (b) at (-0.3,-1.0) {};
    \node[treenode] (c) at (0.3,-1.0) {};
    \draw[treeedge] (r) -- (a);
    \draw[treeedge] (a) -- (b);
    \draw[treeedge] (a) -- (c);
\end{tikzpicture}
\end{center}

Tree 95 enters a new two-leg configuration with depth-46 chains:

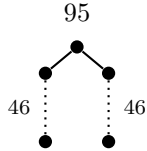
\begin{center}
\begin{tikzpicture}[scale=0.7]
    \node[above] at (0,0.3) {\small 95};
    \node[treenode] (r) at (0,0) {};
    \node[treenode] (l) at (-0.6,-0.5) {};
    \node[treenode] (rt) at (0.6,-0.5) {};
    \draw[treeedge] (r) -- (l);
    \draw[treeedge] (r) -- (rt);
    % Dotted lines for chains
    \draw[thick, dotted] (l) -- (-0.6,-1.8);
    \node[treenode] at (-0.6,-1.8) {};
    \node at (-1.1,-1.15) {\scriptsize 46};
    \draw[thick, dotted] (rt) -- (0.6,-1.8);
    \node[treenode] at (0.6,-1.8) {};
    \node at (1.1,-1.15) {\scriptsize 46};
\end{tikzpicture}
\end{center}

Applying Lemma~\ref{lem:leg} with $x = 46$:
\[
L(46) = 6 \cdot 2^{46} - 2(46) - 6 = 422{,}212{,}465{,}065{,}886
\]

Adding to current step: $95 + 422{,}212{,}465{,}065{,}886 = 422{,}212{,}465{,}065{,}981$.

\begin{center}
\begin{tikzpicture}[scale=0.7]
    \node[above] at (0,0.3) {\scriptsize $422{,}212{,}465{,}065{,}981$};
    \node[treenode] (r) at (0,0) {};
    \node[treenode] (l) at (-0.3,-0.5) {};
    \node[treenode] (rt) at (0.3,-0.5) {};
    \draw[treeedge] (r) -- (l);
    \draw[treeedge] (r) -- (rt);
\end{tikzpicture}
\end{center}

The next tree is a descending chain:

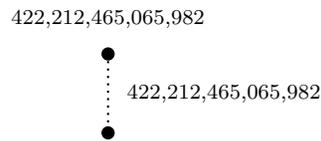
\begin{center}
\begin{tikzpicture}[scale=0.7]
    \node[above] at (0,0.3) {\scriptsize $422{,}212{,}465{,}065{,}982$};
    \node[treenode] (r) at (0,0) {};
    \draw[thick, dotted] (r) -- (0,-1.5);
    \node[treenode] at (0,-1.5) {};
    \node at (2.2,-0.75) {\scriptsize $422{,}212{,}465{,}065{,}982$};
\end{tikzpicture}
\end{center}

\subsection{Final Calculation}

A chain of $n$ nodes can be eliminated in $n - 1$ steps. Therefore:
\[
422{,}212{,}465{,}065{,}982 + (422{,}212{,}465{,}065{,}982 - 1) = 844{,}424{,}930{,}131{,}963
\]

\begin{center}
\begin{tikzpicture}[scale=0.7]
    \node[above] at (0,0.3) {\scriptsize $844{,}424{,}930{,}131{,}963$};
    \node[treenode] (r) at (0,0) {};
\end{tikzpicture}
\end{center}

Since our sequence started at tree 4, we subtract 3:
\[
844{,}424{,}930{,}131{,}963 - 3 = 844{,}424{,}930{,}131{,}960
\]

This completes the proof of Theorem~\ref{thm:main}. \qed

\subsection{Visual Proof of Lemma~\ref{lem:leg}}

This section provides a visual derivation of the Leg Elimination Formula. We show the two-leg elimination process step by step, demonstrating how symmetric configurations reduce through asymmetric intermediate states.

The diagrams below show simplified two-leg trees (without the root and stem) to focus on the elimination pattern. Each pair of vertical chains represents the two legs, and the numbers indicate the tree index in our sequence.

\bigskip

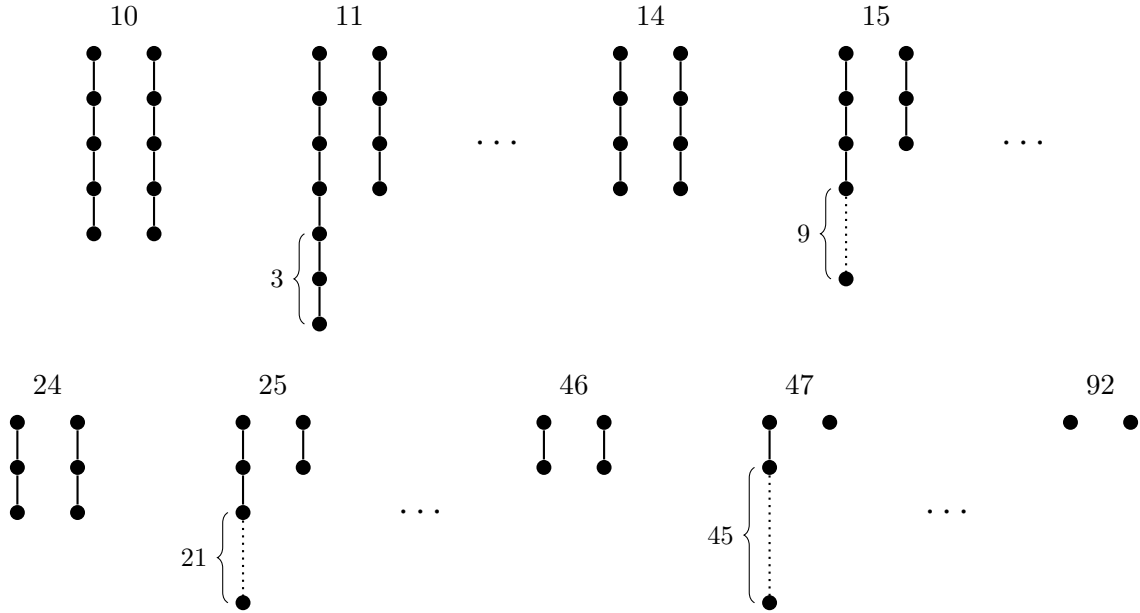
\begin{center}
\begin{tikzpicture}[
    node/.style={circle, fill=black, inner sep=2pt},
    edge/.style={thick}
]

% Label
\node at (0, 0.5) {10};

% Left chain - 5 nodes
\node[node] (l1) at (-0.4, 0) {};
\node[node] (l2) at (-0.4, -0.6) {};
\node[node] (l3) at (-0.4, -1.2) {};
\node[node] (l4) at (-0.4, -1.8) {};
\node[node] (l5) at (-0.4, -2.4) {};
\draw[edge] (l1) -- (l2);
\draw[edge] (l2) -- (l3);
\draw[edge] (l3) -- (l4);
\draw[edge] (l4) -- (l5);

% Right chain - 5 nodes
\node[node] (r1) at (0.4, 0) {};
\node[node] (r2) at (0.4, -0.6) {};
\node[node] (r3) at (0.4, -1.2) {};
\node[node] (r4) at (0.4, -1.8) {};
\node[node] (r5) at (0.4, -2.4) {};
\draw[edge] (r1) -- (r2);
\draw[edge] (r2) -- (r3);
\draw[edge] (r3) -- (r4);
\draw[edge] (r4) -- (r5);

% Tree 11 - shifted to the right
% Label
\node at (3, 0.5) {11};

% Left chain - 7 nodes
\node[node] (l1b) at (2.6, 0) {};
\node[node] (l2b) at (2.6, -0.6) {};
\node[node] (l3b) at (2.6, -1.2) {};
\node[node] (l4b) at (2.6, -1.8) {};
\node[node] (l5b) at (2.6, -2.4) {};
\node[node] (l6b) at (2.6, -3.0) {};
\node[node] (l7b) at (2.6, -3.6) {};
\draw[edge] (l1b) -- (l2b);
\draw[edge] (l2b) -- (l3b);
\draw[edge] (l3b) -- (l4b);
\draw[edge] (l4b) -- (l5b);
\draw[edge] (l5b) -- (l6b);
\draw[edge] (l6b) -- (l7b);
% Brace for bottom 3 nodes
\draw[decorate, decoration={brace, amplitude=4pt, mirror}] (2.4, -2.4) -- (2.4, -3.6) node[midway, left, xshift=-4pt] {\small 3};

% Right chain - 4 nodes
\node[node] (r1b) at (3.4, 0) {};
\node[node] (r2b) at (3.4, -0.6) {};
\node[node] (r3b) at (3.4, -1.2) {};
\node[node] (r4b) at (3.4, -1.8) {};
\draw[edge] (r1b) -- (r2b);
\draw[edge] (r2b) -- (r3b);
\draw[edge] (r3b) -- (r4b);

% Ellipsis
\node at (5, -1.2) {\Large $\cdots$};

% Tree 14 - symmetric 4,4
% Label
\node at (7, 0.5) {14};

% Left chain - 4 nodes
\node[node] (l1c) at (6.6, 0) {};
\node[node] (l2c) at (6.6, -0.6) {};
\node[node] (l3c) at (6.6, -1.2) {};
\node[node] (l4c) at (6.6, -1.8) {};
\draw[edge] (l1c) -- (l2c);
\draw[edge] (l2c) -- (l3c);
\draw[edge] (l3c) -- (l4c);

% Right chain - 4 nodes
\node[node] (r1c) at (7.4, 0) {};
\node[node] (r2c) at (7.4, -0.6) {};
\node[node] (r3c) at (7.4, -1.2) {};
\node[node] (r4c) at (7.4, -1.8) {};
\draw[edge] (r1c) -- (r2c);
\draw[edge] (r2c) -- (r3c);
\draw[edge] (r3c) -- (r4c);

% Tree 15
% Label
\node at (10, 0.5) {15};

% Left chain - 12 nodes shown as 4 + ellipsis + 9
\node[node] (l1d) at (9.6, 0) {};
\node[node] (l2d) at (9.6, -0.6) {};
\node[node] (l3d) at (9.6, -1.2) {};
\node[node] (l4d) at (9.6, -1.8) {};
\draw[edge] (l1d) -- (l2d);
\draw[edge] (l2d) -- (l3d);
\draw[edge] (l3d) -- (l4d);
% Dotted line with brace
\draw[thick, dotted] (l4d) -- (9.6, -3.0);
\node[node] (l5d) at (9.6, -3.0) {};
\draw[decorate, decoration={brace, amplitude=4pt, mirror}] (9.4, -1.8) -- (9.4, -3.0) node[midway, left, xshift=-4pt] {\small 9};

% Right chain - 3 nodes
\node[node] (r1d) at (10.4, 0) {};
\node[node] (r2d) at (10.4, -0.6) {};
\node[node] (r3d) at (10.4, -1.2) {};
\draw[edge] (r1d) -- (r2d);
\draw[edge] (r2d) -- (r3d);

% Ellipsis
\node at (12, -1.2) {\Large $\cdots$};

\end{tikzpicture}
\end{center}

\begin{center}
\begin{tikzpicture}[
    node/.style={circle, fill=black, inner sep=2pt},
    edge/.style={thick}
]

% Tree 24 - symmetric 3,3
% Label
\node at (0, 0.5) {24};

% Left chain - 3 nodes
\node[node] (l1e) at (-0.4, 0) {};
\node[node] (l2e) at (-0.4, -0.6) {};
\node[node] (l3e) at (-0.4, -1.2) {};
\draw[edge] (l1e) -- (l2e);
\draw[edge] (l2e) -- (l3e);

% Right chain - 3 nodes
\node[node] (r1e) at (0.4, 0) {};
\node[node] (r2e) at (0.4, -0.6) {};
\node[node] (r3e) at (0.4, -1.2) {};
\draw[edge] (r1e) -- (r2e);
\draw[edge] (r2e) -- (r3e);

% Tree 25
% Label
\node at (3, 0.5) {25};

% Left chain - 3 visible nodes + 21 shown with dotted line
\node[node] (l1f) at (2.6, 0) {};
\node[node] (l2f) at (2.6, -0.6) {};
\node[node] (l3f) at (2.6, -1.2) {};
\draw[edge] (l1f) -- (l2f);
\draw[edge] (l2f) -- (l3f);
% Dotted line with brace for 21 nodes
\draw[thick, dotted] (l3f) -- (2.6, -2.4);
\node[node] (l4f) at (2.6, -2.4) {};
\draw[decorate, decoration={brace, amplitude=4pt, mirror}] (2.4, -1.2) -- (2.4, -2.4) node[midway, left, xshift=-4pt] {\small 21};

% Right chain - 2 nodes
\node[node] (r1f) at (3.4, 0) {};
\node[node] (r2f) at (3.4, -0.6) {};
\draw[edge] (r1f) -- (r2f);

% Ellipsis
\node at (5, -1.2) {\Large $\cdots$};

% Tree 46 - symmetric 2,2
% Label
\node at (7, 0.5) {46};

% Left chain - 2 nodes
\node[node] (l1g) at (6.6, 0) {};
\node[node] (l2g) at (6.6, -0.6) {};
\draw[edge] (l1g) -- (l2g);

% Right chain - 2 nodes
\node[node] (r1g) at (7.4, 0) {};
\node[node] (r2g) at (7.4, -0.6) {};
\draw[edge] (r1g) -- (r2g);

% Tree 47
% Label
\node at (10, 0.5) {47};

% Left chain - 2 visible nodes + 45 shown with dotted line
\node[node] (l1h) at (9.6, 0) {};
\node[node] (l2h) at (9.6, -0.6) {};
\draw[edge] (l1h) -- (l2h);
% Dotted line with brace for 45 nodes
\draw[thick, dotted] (l2h) -- (9.6, -2.4);
\node[node] (l3h) at (9.6, -2.4) {};
\draw[decorate, decoration={brace, amplitude=4pt, mirror}] (9.4, -0.6) -- (9.4, -2.4) node[midway, left, xshift=-4pt] {\small 45};

% Right chain - 1 node
\node[node] (r1h) at (10.4, 0) {};

% Ellipsis
\node at (12, -1.2) {\Large $\cdots$};

% Tree 92 - symmetric 1,1
% Label
\node at (14, 0.5) {92};

% Left chain - 1 node
\node[node] (l1i) at (13.6, 0) {};

% Right chain - 1 node
\node[node] (r1i) at (14.4, 0) {};

\end{tikzpicture}
\end{center}

\vspace{0.5cm}

\noindent\textbf{Pattern Analysis (with $n = 3$):}

\medskip

\noindent The extra nodes added at each asymmetric step follow the pattern:
\begin{align*}
3 &= n \\
9 &= 2n + n \\
21 &= 4n + 2n + n \\
45 &= 8n + 4n + 2n + n
\end{align*}

\noindent This is the sum $\displaystyle\sum_{i=0}^{m} 2^i \cdot n = (2^{m+1} - 1) \cdot n$

\medskip

\noindent The number of steps between symmetric configurations:
\begin{align*}
10 \to 14: \quad 4 &= n + 1 \\
14 \to 24: \quad 10 &= 2n + n + 1 \\
24 \to 46: \quad 22 &= 4n + 2n + n + 1 \\
46 \to 92: \quad 46 &= 8n + 4n + 2n + n + 1
\end{align*}

\noindent This is the sum $\displaystyle\sum_{i=0}^{m} 2^i \cdot n + 1 = (2^{m+1} - 1) \cdot n + 1$

\bigskip

\noindent\textbf{Deriving the Leg Elimination Formula:}

\medskip

\noindent Substituting $n = 3$ and summing all steps from $m = 0$ to $m = x-1$:
\begin{align*}
\sum_{m=0}^{x-1} \left[(2^{m+1} - 1) \cdot 3 + 1\right] &= 4 + 10 + 22 + 46 + \cdots \\[6pt]
&= \sum_{m=0}^{x-1} \left(3 \cdot 2^{m+1} - 2\right) \\[6pt]
&= 3 \sum_{m=0}^{x-1} 2^{m+1} - 2x \\[6pt]
&= 3 (2^{x+1} - 2) - 2x \\[6pt]
&= 6 \cdot 2^x - 2x - 6
\end{align*}

\noindent where we used the geometric series $\displaystyle\sum_{m=0}^{x-1} 2^{m+1} = 2 + 4 + \cdots + 2^x = 2^{x+1} - 2$.

\medskip

\noindent This completes the derivation of the Leg Elimination Formula: $L(x) = 6 \cdot 2^x - 2x - 6$. \qed

\subsection{A Note on the Initial Embedding}

A point of contention has arisen regarding whether the first tree in our sequence (a root with 3 leaf children) inf-embeds into the second tree (a root with 2 children, where the right child has 2 leaf children). If such an embedding existed, the sequence would be invalid from the start. I clarify here why no such embedding exists.

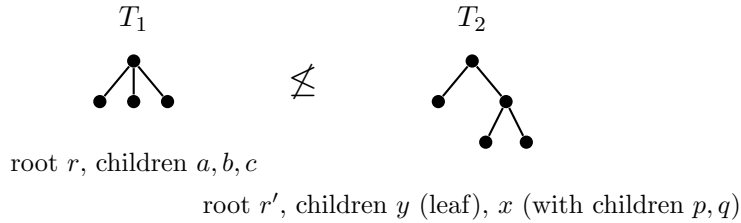
\begin{center}
\begin{tikzpicture}[scale=0.9]
    % Tree 1
    \begin{scope}[shift={(0,0)}]
        \node[above] at (0,0.3) {$T_1$};
        \node[treenode] (r1) at (0,0) {};
        \node[treenode] (a1) at (-0.5,-0.6) {};
        \node[treenode] (b1) at (0,-0.6) {};
        \node[treenode] (c1) at (0.5,-0.6) {};
        \draw[treeedge] (r1) -- (a1);
        \draw[treeedge] (r1) -- (b1);
        \draw[treeedge] (r1) -- (c1);
        \node[below] at (0,-1.2) {\small root $r$, children $a,b,c$};
    \end{scope}

    % Arrow
    \node at (2.5,-0.3) {\Large $\not\leq$};

    % Tree 2
    \begin{scope}[shift={(5,0)}]
        \node[above] at (0,0.3) {$T_2$};
        \node[treenode] (r2) at (0,0) {};
        \node[treenode] (y2) at (-0.5,-0.6) {};
        \node[treenode] (x2) at (0.5,-0.6) {};
        \node[treenode] (p2) at (0.2,-1.2) {};
        \node[treenode] (q2) at (0.8,-1.2) {};
        \draw[treeedge] (r2) -- (y2);
        \draw[treeedge] (r2) -- (x2);
        \draw[treeedge] (x2) -- (p2);
        \draw[treeedge] (x2) -- (q2);
        \node[below] at (0,-1.8) {\small root $r'$, children $y$ (leaf), $x$ (with children $p,q$)};
    \end{scope}
\end{tikzpicture}
\end{center}

\begin{proof}[Proof that $T_1 \not\leq T_2$]
Suppose for contradiction that an inf-embedding $f: V(T_1) \to V(T_2)$ exists.

In $T_1$, note that $\inf(a, b) = \inf(b, c) = \inf(a, c) = r$, since all three children share only the root as their common ancestor.

\textbf{Case 1:} $f(r) = r'$. Then $a, b, c$ must map to descendants of $r'$. The available leaves in $T_2$ are $\{y, p, q\}$---exactly three. So we must have $\{f(a), f(b), f(c)\} = \{y, p, q\}$ in some order.

Without loss of generality, suppose $f(a) = p$ and $f(b) = q$. Then:
\[
f(\inf(a, b)) = f(r) = r'
\]
But the inf-embedding condition requires:
\[
f(\inf(a, b)) = \inf(f(a), f(b)) = \inf(p, q) = x \neq r'
\]
This is a contradiction. Any other assignment of $\{a,b,c\}$ to $\{y,p,q\}$ similarly fails, since at least two of them will be mapped to $\{p, q\}$, whose infimum is $x$, not $r'$.

\textbf{Case 2:} $f(r) = x$. Then $a, b, c$ must all map to descendants of $x$. But $x$ has only two descendants ($p$ and $q$), so we cannot injectively map three vertices.

\textbf{Case 3:} $f(r) \in \{y, p, q\}$. These are leaves with no descendants, so we cannot map the children $a, b, c$ of $r$ to descendants.

All cases lead to contradiction. Therefore $T_1 \not\leq T_2$.
\end{proof}

The key insight is that in $T_1$, any two of the three children share the root as their infimum, but in $T_2$, two of the three available leaves ($p$ and $q$) share $x$ as their infimum rather than the root. This structural difference prevents any valid inf-embedding.

\section{Discussion}

\subsection{Comparison with Previous Estimates}

Our bound of $\mathrm{tree}(3) \geq 3 \cdot 2^{48} - 8 \approx 8.44 \times 10^{14}$ significantly exceeds Friedman's conjecture that $\mathrm{tree}(3) < 100$ \cite{friedman2006}. This discrepancy likely arises from the non-obvious nature of the optimal opening sequence. Naive constructions easily yield sequences of length 100--200, and discovering the configuration that triggers exponential growth requires careful analysis of the leg elimination process.

\subsection{Upper Bounds and Exactness}

It remains an open question whether this lower bound is tight. I conjecture:

\begin{conjecture}
$\mathrm{tree}(3) = 844{,}424{,}930{,}131{,}960$.
\end{conjecture}

Proving this would require showing that no alternative initial sequence can yield a longer valid sequence---a challenging combinatorial task that likely requires exhaustive case analysis of all possible opening moves.

\subsection{Implications for tree(4) and Beyond}

Kihara \cite{kihara2020} has shown that $\mathrm{tree}(4) > G$, where $G$ is Graham's number, and more precisely that $\mathrm{tree}(4) > f_{\varepsilon_0}(G)$ in the fast-growing hierarchy. The dramatic jump from $\mathrm{tree}(3) \approx 10^{14}$ to $\mathrm{tree}(4) > f_{\varepsilon_0}(G)$ illustrates the extraordinary growth rate of the weak tree function. Even the relatively ``tame'' value of $\mathrm{tree}(3)$ already exceeds what one might naively expect, hinting at the explosive behavior that emerges for larger arguments.

\section*{Acknowledgements}

I thank the Stack Exchange community for helpful discussions on combinatorics and tree embeddings. I am grateful to Professor Takayuki Kihara for his foundational work on lower bounds for the weak tree function, and for his suggestions on improving this paper, including the use of tree diagrams for better visualization and the addition of the proof of Lemma~\ref{lem:leg}.

\end{document}